\documentclass[11pt]{article}

\usepackage{amsfonts}
\usepackage{epsfig,graphics,subfigure,psfrag,amsmath,amssymb}
\usepackage{latexsym,wrapfig}
\usepackage{bm}
\usepackage{amsthm}
\usepackage{graphicx}
\usepackage{mathrsfs}
\usepackage{fancyref}
\usepackage{dcolumn}
\usepackage{color}
\usepackage{pstricks}
\usepackage{epstopdf}
\usepackage{subfigure}
\usepackage{multicol}
\usepackage{enumerate}
\usepackage{caption}
\usepackage{geometry}
\usepackage{threeparttable}

\usepackage{hyperref}
\makeatletter
\@addtoreset{equation}{section}
\makeatother

\usepackage{latexsym,amsmath,amssymb,paralist}

\usepackage{comment}
\usepackage[compress]{cite}
\usepackage{algorithm}
\usepackage{algpseudocode}
\usepackage[title]{appendix}

\usepackage{graphicx}
\usepackage{epstopdf}
\usepackage{bm}
\usepackage{hyperref}
\usepackage[normalem]{ulem}

\newtheorem{thm}{Theorem}[section]

\newtheorem{defn}[thm]{Definition}

\newtheorem{rem}[thm]{Remark}
\def\E\mathbb{ E}

\newcommand{\be}{\begin{eqnarray}}
\newcommand{\ee}{\end{eqnarray}}
\newcommand{\ben}{\begin{eqnarray*}}
\newcommand{\een}{\end{eqnarray*}}

\def\XXint#1#2#3{{\setbox0=\hbox{$#1{#2#3}{\int}$}
\vcenter{\hbox{$#2#3$}}\kern-.51\wd0}}

{\allowdisplaybreaks
\begin{document}
\title{The Cumulative Distribution Function Based Method for Random Drift Model}
\author{Chenghua Duan \footnotemark[1] \footnotemark[2]\and Chun Liu\footnotemark[3]\and Xingye Yue \footnotemark[4] }
%
 \renewcommand{\thefootnote}{\fnsymbol{footnote}}
 \footnotetext[1]{Department of Mathematics, Shanghai University, Shanghai, China. Email:  chduan@shu.edu.cn.}

 \footnotetext[2] {Newtouch Center for Mathematics of Shanghai University, Shanghai, China, 20044.}

  \footnotetext[3]{Department of Applied Mathematics, Illinois Institute of Technology, Chicago, IL 60616, USA. Email: cliu124@iit.edu.}
  \footnotetext[4]{School of Mathematical Sciences, Soochow University, Suzhou 215006, Jiangsu, China. Email: xyyue@suda.edu.cn.}

 %
\maketitle

\begin{abstract}
In this paper, we propose a numerical method   to uniformly handle the random genetic drift model for pure drift with or without natural selection and mutation.
 For pure drift and natural selection case, the Dirac $\delta$ singularity will develop at two boundary ends and the mass lumped at the two ends stands for the fixation probability. For the one-way mutation case, known as  Muller's ratchet,  the accumulation of deleterious mutations leads to the loss of the fittest gene,   the Dirac $\delta$ singularity  will spike only at one  boundary end, which stands for the fixation of the deleterious gene and loss of the fittest one. For  two-way mutation case,  the singularity with negative power law may emerge near boundary points.   We first rewrite the original model on the probability density function (PDF) to one with respect to the cumulative distribution function (CDF). Dirac  $\delta$ singularity of the  PDF becomes the discontinuity of the CDF.  Then we establish a  upwind scheme, which  keeps the  total probability, is positivity preserving and  unconditionally stable. For pure drift, the scheme  also  keeps the conservation of expectation.  It can catch the discontinuous jump of the CDF, then predicts accurately  the fixation probability for pure drift with or without natural selection and one-way mutation. For two-way mutation case, it can catch the power law of the singularity.  
The numerical results show  the effectiveness  of the scheme.

\bigskip


\bigskip
\noindent

\end{abstract}
\section{Introduction}
In the population genetics, the random genetic drift model describes that  the number of gene variants (alleles) fluctuates randomly over time  due to random sampling.  The fraction of an allele  in the population can be used to measure the intensity of the random genetic drift.  In other words,  the value of the  fraction equals to zero or one, which means the allele disappears or is completely  chosen in the system.
Population genetics models, aiming at modeling genetic variability, had a natural start with discrete stochastic models at the individual
level. The earlier mathematical model of random genetic drift,  known as the Wright-Fisher model,  was constructed by Ronald Fisher \cite{R.A.Fisher(1923)} and Sewall Wright \cite{S.Wright(1929), S.Wright(1937),  S.Wright(1945)}.  The Wright-Fisher model is regarded as a discrete-time Markov chain under the assumption that the generations do not overlap and that each copy
of the gene of the new generation is selected independently and randomly from the whole gene pool of the previous generation.  Then Moran \cite{Moran(1958),Moran(1962)} and Kimura \cite{M.Kimura(1955a),M.Kimura(1955b),M.Kimura(1962),M.Kimura(1964)} derived the diffusion limit of random Wright-Fisher model.
Chalub et. al spreaded the large population limit of the Moran process, and obtained a continuous model that may highlight the genetic-drift (neutral evolution) or natural selection \cite{M.Kimura(1964),Chalub(2009),Traulsen(2005)}.

We consider a population of constant size $N_e$ with a pair of two types A and B.  Time $t$
is increased by the time step $\Delta t$, and the process is repeated. $x$ is the gene $A$ frequency at $t$ generation. 
 Let $P_{N_e,\Delta t}(t,x)$ denote the probability
of finding a fraction $x$ of type $A$ individuals at time $t$ in population
of fixed size $N_e$ evolving in discrete time, according to the Moran
process. Then, in the limit of large
population and small time steps, we postulate the existence of a
probability density of states, that will depend on the precise way
the limits are taken \cite{Chalub(2009)}. Namely:
$$f(t,x)=\lim\limits_{N_e\rightarrow\infty,\\ \Delta t\rightarrow 0}\frac{P_{N_e,\Delta t}(t,x)}{1/N_e}.$$


In this paper, we study the following initial-boundary problem of $f(t,x)$:
\begin{equation}\label{equ:forwardK}
\left\{
 \begin{aligned}
 & \partial_t f-\partial_{xx}[x(1-x)f]+\partial_x[M(x)f]=0,\ x\in\Omega=(0,1),\ t>0,\\
 & \mathcal{J}=-\partial_x (x(1-x)f)+M(x)f, \ \mathcal{J}\cdot{\bf n}|_{\partial \Omega}=0,\\
 &f(0,x)=f_0(x), x\in\Omega,
 \end{aligned}
 \right.\
\end{equation}
where the function $M(x)$, which is typically a polynomial in $x$, incorporates the
forces of migration, mutation, and selection,  acting at time t\cite{Zhao(2013),M.Kimura(1964),M.Kimura(1955a)}.


Eq.\eqref{equ:forwardK} can be supplemented by the following conservation laws:
\begin{itemize}
  \item[(a)] Mass conservation law:
  \begin{equation}
\label{equ:massConservation}
\int_0^1f(0,x)dt=\int_0^1f(t,x)dt=1.
\end{equation}
  \item[(b)]  Moment conservation law:
  \begin{equation}
\label{equ:Conservation2}
\frac{d}{dt}\int_0^1\theta(x)f(t,x)dx=0,
\end{equation}
where $\theta(x)$  is the fixation probability function that satisfies
\begin{equation}\label{equ:theta}
x(1-x)\theta''+M(x)\theta'=0, \ \theta(0)=0,\ \theta(1)=1.\end{equation}
\end{itemize}

One notes that    \eqref{equ:Conservation2}  recovers the conservation of   expectation for pure drift ($M(x)=0$) with $\theta(x)=x$.

Based on   different   $M(x)$, the behavior of the solution  has three cases:

 \noindent{\bf Case 1. Pure drift and natural selection.\ }
For pure drift, we have $M(x)=0$. When natural selection is involved, we have \cite{Chalub(2009a),Carrillo(2022)},
\begin{equation}\label{equ:MS}
M(x)=x(1-x)(\eta x +\beta).
\end{equation}
 The important theoretically relevant results  were shown in \cite{Chalub(2009),Chalub(2009a),Carrillo(2022)}. Let $\mathcal{BM}^+([0,1])$ denote
the space of (positive) Radon measures on $[0,1]$.
\begin{defn}\label{def:weakSolution}
A weak solution to  \eqref{equ:forwardK} is a function in $L^{\infty}([0,\infty),\mathcal{BM}^+([0,1]))$   satisfying
\begin{equation*}
\begin{split}
-\int_0^{\infty}\int_0^1f(t,x)\partial_t\phi(t,x)dxdt
&=\int_0^{\infty}\int_0^1f(t,x)x(1-x)\left(\partial_{xx} \phi(t,x)+\frac{M(x)}{x(x-1)}\partial_x \phi(t,x)\right)dxdt\\
& +\int_0^1f_0(x)\phi(0,x)dx,\  \mbox{for\ any\ } \phi(t,x)\in C_c^{\infty}([0,+\infty)\times[0,1]).
\end{split}
\end{equation*}
\end{defn}
 The following theorem can be found in \cite{Chalub(2009a)}.
\begin{thm}\label{thm}
 If $f_0\in\mathcal{BM}^+([0,1])$ and $M(x)=0$ or given by \eqref{equ:MS}, then there exists a unique weak solution $f\in L^{\infty}([0,\infty);\mathcal{BM}^+([0,1]))\cap C^{\infty}(\mathbb{R}^+,C^{\infty}((0,1)))$
to equation \eqref{equ:forwardK}, in a sense of Definition \eqref{def:weakSolution}, such that  the conservation laws \eqref{equ:massConservation}-\eqref{equ:Conservation2} are valid. The solution has a form as
 $$f(t,x)=r(t,x)+a(t)\delta_0+b(t)\delta_1,$$
where $\delta_s$ denotes the singular point measure supported at $s$,  $r(t,x)\in C^{\infty}(\mathbb{R}^+;C^{\infty}([0,1]))  $ is
a classical solution to \eqref{equ:forwardK} without any boundary condition, $a(t)$ and $b(t)$   $\in C([0,\infty))\cap C^{\infty}(\mathbb{R}^+)$.   
 Furthermore, $\lim\limits_{t\rightarrow\infty}r(t,x)=0$ uniformly,  and $a(t)$ and $b(t)$ are monotonically increasing functions such that
$$a^{\infty}:=\lim\limits_{t\rightarrow\infty}a(t)=\left(1-\int_0^1f_0(x)\theta(x)dx\right),$$
and
$$b^{\infty}:=\lim\limits_{t\rightarrow\infty}b(t)=\int_0^1f_0(x)\theta(x)dx,$$
where $\theta(x)$  is given by \eqref{equ:theta}.  Finally, the convergence rate is exponential.
\end{thm}

The appearance of the point measure $\delta_0$ ($\delta_1$) stands for that the fixation at gene $B$ ($A$) happens with a probability $a(t)$ ($b(t)$).
The spatial temporal dynamics of the Kimura equation are well understood in the purely diffusive case and in only a relatively small number of population
\cite{M.Kimura(1955),Crow(1970),Crow(1956)}.

 \noindent{\bf Case 2.   One-way mutation: Muller's ratchet.}\ Considering the one-way mutation from gene $B$ to gene $A$,
\begin{equation}\label{equ:MO}
M(x)=\gamma(1-x),
\end{equation}
where the constant  $\gamma>0$ stands for the mutation rate \cite{M.Kimura(1955a),Ewens(2004),Uecker(2011),Kaushik(2023)}.

A well-known model is Muller's ratchet, i.e, the fittest gene $B$ of individuals is eventually lost from the population and deleterious mutations ($B \rightarrow A$) slowly but irreversibly accumulate through rare stochastic
fluctuations  \cite{Muller(1964),Felsenstein(1974)}.  In a finite asexual population,  offspring inherit all the deleterious mutations their parents possess. Since these offspring also occasionally acquire new deleterious mutations, populations will
tend to accumulate deleterious mutations over time. This effect is known as Muller's ratchet.

There exists a  unique steady solution $f_{\infty}$,
\begin{equation}
\label{equ:SteadyMO}
  f_{\infty}=\delta_1, \ i.e.,\ \mbox{the\ fittest\ gene\ is\ complitely\ lost with probability $1$}.
 \end{equation}

 \noindent{\bf Case 3.  Two-way mutation.}
\begin{equation}\label{equ:MT}
   M(x)=\gamma(1-x)-\mu x,\ \gamma,\mu\in(0,1),
\end{equation}
where the constant  $\mu>0$ stands for the mutation rate of gene $A$ to gene $B$ and $\gamma>0$ is the rate in opposite direction.
In the long term   one might expect   to exist  an equilibrium  state due to the  two direction mutation.  Actually, there exists a  unique steady state solution $f_{\infty}(x)$ to \eqref{equ:forwardK} and \eqref{equ:MT},
\begin{equation}\label{MutationS}f_{\infty}(x)=\frac{C}{x^{1-\gamma}(1-x)^{1-\mu}},\ x\in(0,1),\end{equation}
where $C$ is the   constant such that $\int_0^1f_{\infty}(x)dx=1$. Only singularity of negative power law  develops at the ends, rather than Dirac $\delta$ appears for the cases of pure drift, natural selection and one-way mutation.


For the  numerical simulation,  the crucial
features  to solve \eqref{equ:forwardK} are that the numerical solution can keep the total probability conservation law  and accurately capture the concentration phenomena at the discrete level.   The total probability fails to keep the total probability by some classical numerical schemes \cite{Barakat(1978),Wang(2004),M.Kimura(1955a)}. A complete solution, whose total probability is unity, obtained
by
finite volume method (FVM) schemes in \cite{Zhao(2013)}.  Xu et al. \cite{Xu(2019a)}
compared  a serial of finite volume and finite element schemes for the pure diffusion equation. Their critical comparison of the long-time asymptotic performance urges carefulness in choosing a numerical method for this type of problem, otherwise the
main properties of the model might be not kept.  In recent years,     a variational particle method was proposed based on an energetic variational approach, by which a complete numerical solution can be obtained and  the positivity of the solution  can be kept \cite{Duan(2019)}. However, some artificial criteria must be introduced to handle the concentrations near the boundary ends, even though it is designed based on the biological background.
 Recently, an optimal mass transport method  based on  pseudo-inverse of CDF is used to solve the model \cite{Carrillo(2022)}.   In this method,
the feasible solution is strictly monotonous. However, for the cases of pure drift, natural selection and one-way mutation, the Dirac-$\delta$ concentrations must be developed, then  the corresponding CDF  must be discontinuous at the concentration points. This leads to a fact that its pseudo-inverse can not be strictly monotonous. However, numerical results were presented for the cases of pure drift and selection. So  some manual intervention must be introduced in the numerical implements.    The pure drift problem with multi-alleles, corresponding to a multidimensional PDE, is
investigated in  \cite{Xu(2019b)} by finite-difference methods, where the authors proposed a
numerical scheme with absolute stability and several biologically/physically
relevant quantities  conserved, such as positivity, total probability, and expectation. So far, although quite a few  numerical methods have been established for the pure drift and selection cases, efficient numerical  works on the mutation case  are not reported yet.
%

In this paper, we rewrite the system \eqref{equ:forwardK} to the following new one on the cumulative distribution function (CDF) $F(t,x)=\int_0^x f(t,s)ds$ as
\begin{equation}\label{equCDF}
\left\{
 \begin{aligned}
 & \partial_t F-\partial_x[x(1-x)\partial_x F]+M(x)\partial_x F=0,\ x\in\Omega=(0,1),\ t>0,\\
 & F(t,0^-)=0,\ F(t,1^+)=1,\\
 &F(0,x)=\int_0^x f_0(s)ds, x\in\Omega.
 \end{aligned}
 \right.\
\end{equation}
Taking the CDF into account, the Dirac $\delta$ singularity at the boundary points of original PDF $f$ will change to the discontinuity of the CDF $F$ at the boundary points,  and the fixation probability (lumped mass) will change to the height of the discontinuous jump. The singularity of negative power law will change to a bounded positive power law.

  We will propose a  upwind  numerical scheme for \eqref{equCDF} with a key revision near the boundary, which can  handle the pure drift with or without selection and mutation, is unconditional stable, and keeps  the total probability and positivity.  It also keeps  the conservation of expectation for pure drift.    The numerical results show that the  scheme can catch the height of discontinuity at the ends and predict accurately the fixation probability for the cases of pure drift, natural selection and one-way mutation.  It also predict accurately the negative power of the power law for two-way mutation case.

 The rest of this paper is organized as follows. In Section \ref{sec:num}, we construct the numerical scheme for  \eqref{equCDF}. Some numerical analysis is presented in Section \ref{sec:ana}. In Section \ref{sec:numericalResults}, several numerical examples  are presented to validate the theoretical results and to demonstrate  the ability to trace
the long-time dynamics of random genetic drift. Some discussions will be presented in Section \ref{sec:FDM2} about the relations between the revised scheme and the standard  finite difference method.
\section{Numerical Scheme}\label{sec:num}
In this section, we introduce  the revised finite difference method (rFDM)  for \eqref{equCDF} with the  central difference method for diffusion term and the upwind scheme for convection term.   Let  $h=1/K$  be the spatial step size, and $x_i=ih$, $i=0,\cdots,{K}$, be the spatial  grid points.  Let $\tau$ be the temporal step size, and  $t_n=n\tau$, $n=1,2,\cdots,N$, be the temporal  grid points.

 Define the first order difference as
 \begin{equation}\label{D}
 D_h F_i=\frac{F_{i}-F_{i-1}}{h}, \ \ D_h^{up} F_i=
\left\{
 \begin{aligned}
 &D_h F_i,\  M(x_i)>0 ,\\
 &D_h F_{i+1},\ M(x_i)<0,
 \end{aligned}
 \right.\,\end{equation}
and denote the diffusion coefficient by $a(x)=x (1-x)$.
The upwind numerical scheme, referred as {\bf rFDM}, reads as: Given $F^{n-1}$,  $F^{n}$  solves the following linear algebra system,
\begin{equation}\label{equ:FDM}
\frac{F^{n}_i-F^{n-1}_i}{\tau}-\frac{a_{i+\frac{1}{2}}D_hF^{n}_{i+1}-a_{i-\frac{1}{2}}D_hF^{n}_{i}}{h}+M(x_i) D^{up}_h F^{n}_i=0,\ \ i=1,\cdots,{K}-1,
\end{equation}
with $F^{n}_0=0$ and $F^{n}_{K}=1$, for $n=1,\cdots,N$, and the key revision on the diffusion cofficient
\begin{equation}\label{equ:aboundary} a_{i-\frac 12}= a(x_{i-\frac 12}), \ i=2,\cdots,{K}-1,\ \mbox{\bf but}\ a_{\frac{1}{2}}=a_{{K}-\frac{1}{2}}=0.
 \end{equation}

Finally, the solution of  the original equation  \eqref{equ:forwardK} can be recovered from $F$ by, for $n=1,\cdots,N$,
\begin{equation}\label{equf}
f_i^n=
\left\{
 \begin{aligned}
 & \frac{F_{i+1}^n-F_{i-1}^n}{2h},\ i=2,\cdots,K-2,\\
 & \frac{F_{i+1}^n-F_i^n}{h},\  i=0, 1,\\
 & \frac{F_{i}^n-F_{i-1}^n}{h}, \ i=K-1,K.
 \end{aligned}
 \right.\
\end{equation}
In the above formula, central difference does not applied at near boundary points $i=1$ or $K-1$, thanks to the fact that the discontinuous jump may occur at the boundary points.
\begin{rem}
For problem \eqref{equCDF} in continuous sense, the diffusion coefficient $a(0)=a(1)=0$ degenerates at the boundary points. This means that the information at the boundary points can never be transferred into the domain $\Omega=(0,1)$ by diffusion. In our revised treatment  \eqref{equ:aboundary} for numerical scheme, we set $a_{\frac{1}{2}}=a_{{K}-\frac{1}{2}}=0$, where  a term of $O(h)$ order is omitted, since the exact value $a(x_{\frac{1}{2}})=a(x_{K-\frac{1}{2}})=\frac h2 (1-\frac h2)$. With this revised treatment,  the boundary value $F_0^n$ and $F_K^n$ is not involved in the discrete system \eqref{equ:FDM} if $M\equiv 0$, i.e., the boundary value can never be transferred into the inner points by discrete diffusion. In   Section \ref{sec:FDM2}, we will discuss the standard scheme without this revision.
\end{rem}

\section{Analysis Results}\label{sec:ana}
 In this section, we will focus on the numerical analysis for rFDM \eqref{equ:FDM}, including the unconditional stability, positivity preserving,  and conservation law of the total probability and expectation.
\begin{thm}
The  upwind scheme rFDM \eqref{equ:FDM} is unconditionally stable and positivity preserving.
\end{thm}
\noindent{\bf Proof:}
Without loss of generality, assume there exists $i^{*}$ such that $M(i)>0$, $i=1,\cdots,i^*$ and $M(i)<0$, $i=i^*+1,\cdots,K-1$.    \eqref{equ:FDM} can be written as:
\begin{equation}\label{equ:linearsystem}
\begin{split}
&\frac{F^{n}_i-F^{n-1}_i}{\tau}
-\frac{a_{i+\frac{1}{2}}F^{n}_{i+1}-(a_{i+\frac{1}{2}}+a_{i-\frac{1}{2}})F^{n}_{i}+a_{i-\frac{1}{2}}F^{n}_{i-1}}
{h^2}+M(x_i)\frac{F^{n}_i-F^{n}_{i-1}}{h}=0,\\&\ i=1,\cdots i^{*},\\
&\frac{F^{n}_i-F^{n-1}_i}{\tau}
-\frac{a_{i+\frac{1}{2}}F^{n}_{i+1}-(a_{i+\frac{1}{2}}+a_{i-\frac{1}{2}})F^{n}_{i}+a_{i-\frac{1}{2}}F^{n}_{i-1}}
{h^2}+M(x_i)\frac{F^{n}_{i+1}-F^{n}_{i}}{h}=0,\\&\ i=i^{*}+1,\cdots,K-1.
\end{split}
\end{equation}
with $F^{n}_{0}=0$, $F^{n}_{K}=1$ at time $t^n$, $n=1,\cdots,N$.

Let ${\bf B}= (b_{ij})$ be the matrix of the linear system. Then ${\bf B}$ is a tri-diagonal matrix.

\noindent For $i=1,\cdots,i^*$,
$$b_{i i}=1+\frac{\tau}{h^2}(a_{i+\frac{1}{2}}+a_{i-\frac{1}{2}})+\frac{\tau}{h}M(x_i),$$
$$b_{i,i-1}=-\frac{\tau}{h^2}a_{i-\frac{1}{2}}-\frac{\tau}{h}M(x_i),$$
$$b_{i,i+1}=-\frac{\tau}{h^2}a_{i+\frac{1}{2}}.$$
For $i=i^*+1,\cdots,K-1$,
$$b_{i i}=1+\frac{\tau}{h^2}(a_{i+\frac{1}{2}}+a_{i-\frac{1}{2}})-\frac{\tau}{h}M(x_i),$$
$$b_{i,i-1}=-\frac{\tau}{h^2}a_{i-\frac{1}{2}},$$
$$b_{i,i+1}=-\frac{\tau}{h^2}a_{i+\frac{1}{2}}+\frac{\tau}{h}M(x_i).$$
For $i=0$,
$b_{i i}=1,\ b_{i,i+1}=0.$
For $i=K$,
$b_{i i}=1, \  b_{i,i-1}=0.$

 Note that  $$b_{i i}>0, \  i=0,\cdots,K$$   $$b_{i,i-1}\leq0, \ i=1,\cdots,K,$$ $$b_{i,i+1}\leq0,\ i=0,\cdots,K-1,$$
and $$b_{ii}-\Big||b_{i,i-1}|+|b_{i,i+1}|\Big|=1>0,\ i=1,\cdots,K-1,$$
Then ${\bf B}$ is a M-matrix, i.e., any entry of the inverse matrix ${\bf B}^{-1}$ is positive. So rFDM \eqref{equ:FDM} is unconditionally stable and positivity preserving.

\begin{thm}
The numerical scheme rFDM \eqref{equ:FDM}  keeps    the conservation  of total probability.  For pure drift case,   the conservation of  expectation  also holds.
\end{thm}
\noindent{\bf Proof:}
Define $p_i^n := F_i^n - F_{i-1}^n$ be the probability  for the fraction of gene $A$ belongs to the interval $I_i=[x_{i-1}, x_i]$. Then the total probability is
$$P_{total}^n= \sum\limits_{i=1}^{K} p^n_i = F^n_K - F^n_0 = 1, \ \ n=0,1,\cdots, N,$$
i.e., the discrete system \eqref{equ:FDM} keeps the mass conservation naturally.

By the integral by parts,  the expectation $\mathcal{E}(t)$ satisfies that
\begin{equation} \label{exp-dis} \mathcal{E}(t):=\int_0^1 xf(t,x)dx=\int_0^1 x\partial_x F(t,x)dx=1-\int_0^1 F(t,x)dx.$$
  So we define a  discrete expectation   as  $$\mathcal{E}^n_h:=1-\sum\limits_{i=1}^{K-1}F^n_ih-\frac{h}{2}F^n_0-\frac{h}{2}F^n_K.\end{equation}

For pure drift case, $M(x)\equiv 0$, then  we have, by \eqref{equ:FDM} and \eqref{equ:aboundary}, that
\begin{equation*}
\begin{aligned}
\mathcal{E}^{n}_h-\mathcal{E}^{n-1}_h
&=-\sum\limits_{i=1}^{K-1}h(F_i^{n}-F_{i}^{n-1})=-\sum\limits_{i=1}^{K-1}\tau\left(a_{i+\frac{1}{2}}D_hF^{n}_{i+1}-a_{i-\frac{1}{2}}D_hF^{n}_{i}\right)\\
&=-\tau(a_{K-\frac{1}{2}}D_hF^{n}_{K}-a_{\frac{1}{2}}D_hF^{n}_{1}) = 0.
\end{aligned}\end{equation*}
So the discrete system \eqref{equ:FDM}   keeps the conservation of the discrete expectation for pure drift case.
Please note that, without the revised treatment \eqref{equ:aboundary}, the conservation of the discrete expectation is invalid.

\section{Numerical results}\label{sec:numericalResults}
In this section, we will show the   effectiveness of this algorithm by different numerical examples. In Example 1, we verify the local convergence.  In Section \ref{sec:numericalResults1},   the  pure drift $M=0$ and natural selection $M(x)=x(1-x)(\eta x+\beta)$   case are studied in  Example 2 and  Example  3, respectively.
In Section \ref{sec:numericalResults3},  we consider the mutation case   $M(x)=\gamma(1-x)-\mu x$, including  one-way mutation model, such as the Muller's ratchet model,  in Example 4 and  two-way mutation model  in  Example 5. 



{\bf Example 1.   Local convergence}

Although the singularity may develop at boundary points, the solution in interior area  is sufficiently smooth. To verify  the correctness of the numerical scheme, we check the local convergence.

Let $\Omega'\subset\Omega$ be the interior area and
$$k_1:=\inf\limits_{0\leq i \leq K}\{i|x_i\in\Omega'\},$$
$$k_2:=\max\limits_{0\leq i \leq K}\{i|x_i\in\Omega'\}.$$
Define
the error $e$ of $F(t,x)$ in   $\mathcal{L}^2(\Omega')$ and $\mathcal{L}^{\infty}(\Omega')$ mode   as
\begin{equation}\label{equ: error2}
\|e\|_2:=\left(\sum\limits_{i=k_1}^{k_2}(F_i^n-F^n_{e,i})^2h\right)^{\frac{1}{2}},
\end{equation}
and
\begin{equation}\label{equ: errorinf}
\|e\|_{\infty}:=\max\limits_{k_1\leq i \leq k_2}\{|F^n_i-F^n_{e,i}|\},
\end{equation}
where $\{F_i^n\}_{i=0}^K$ is the numerical solution of CDF  model \eqref{equCDF}, and $\{F_{e,i}^n\}_{i=0}^K$ is the corresponding  exact solution  at time $t^n$, $n=1,\cdots,N$.

In this example, we take the initial probability density as
$$f_0(x)=1,\ x\in[0,1].$$
Table \ref{table:Correct} shows the error and local  convergence order of  $F(t,x)$ for $M(x)= x(1-x)(-4x+2)$ and  $M(x)=0.2(1-x)+0.4x$ in $\mathcal{L}^2([0.3,0.7])$ and $\mathcal{L}^{\infty}([0.3,0.7])$ mode with different space  and time grid sizes  at time $t=0.1$.    The
reference "exact" solution is obtained numerically on a  fine mesh with   $h=1/100000$ and $\tau=1/100000$. The results show that the local convergence of the numerical scheme is 2nd order   in space and  1st order  in time for different $M(x)$ in the inner region $[0.3,0.7]$.
\begin{table}[ht]
\centering
\footnotesize{\caption{Error of $F$  in $\Omega'=[0.3,0.7]$ at   $t=0.1$  in {\bf Example 1}.}\label{table:Correct}
\begin{tabular}{cccccc}
\hline
 &  &  $M(x)= x(1-x)(-4x+2)$& & &\\
\hline
 $h$    &$\tau$ &$\|e\|_2$ &order   &$\|e\|_{\infty}$&order\\\hline
 1/100&1/100 &9.78093e-04& &2.61509e-03& \\\hline
 1/200 &1/400 &2.41752e-04&2.01644&6.62837e-04&1.98013\\\hline
 1/400 &  1/1600 &6.03394e-05&2.00236&1.69222e-04&1.96973\\

 \hline
 &   & $M(x)=0.2(1-x)+0.4x$& & &\\
\hline

 $h$    &$\tau$ &$\|e\|_2$ &order   &$\|e\|_{\infty}$&order\\\hline
 1/100&1/100 &9.90565e-04& &3.05562e-03& \\\hline
 1/200 &1/400 &2.47424e-04&2.00127&7.76938e-04&1.97559\\\hline
 1/400 &  1/1600 &6.16308e-05&2.00526&1.96225e-04&1.98529\\

 \hline
\end{tabular}}
\end{table}
\subsection{Pure drift and natural selection  }\label{sec:numericalResults1}
In this section, we study pure drift and natural selection case which keeps the conservation law \eqref{equ:massConservation} and \eqref{equ:Conservation2}.

 Firstly, we define the error of fixation probability at the left and right boundary points at large enough  time (near the steady state) $T=\tau N$ as
 \begin{equation}
 \label{errorl}
 e_{left}:=|F_1^N-F_0^N-a^{\infty}|,\end{equation}
  and
\begin{equation} \label{errorr}
e_{right}:=|F_K^N-F_{K-1}^N-b^{\infty}|,\end{equation}
where  $a^{\infty}$ and $b^{\infty}$ are the exact fixation probability at boundary points given in Theorem \ref{thm}.

\noindent{\bf Example 2.  Pure drift }

In this example, we consider pure drift case, i.e.,  $M(x)=0$, with a Gaussian distribution initial function   as
\begin{equation}
\label{equ:Nini}
f(0,x)=\frac{1}{\sqrt{2\pi}\sigma}\exp^{-\frac{(x-x_0)^2}{2\sigma^2}},
\end{equation} with $\sigma=0.01$ and $x_0=0.7$.

\begin{figure}[ht]
\centering
{\includegraphics[scale=0.8]{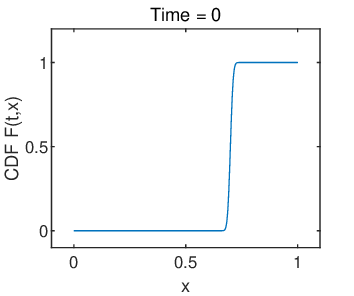}}
{\includegraphics[scale=0.8]{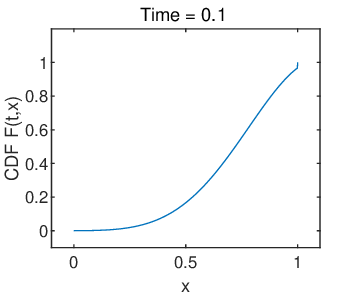}}\\
{\includegraphics[scale=0.8]{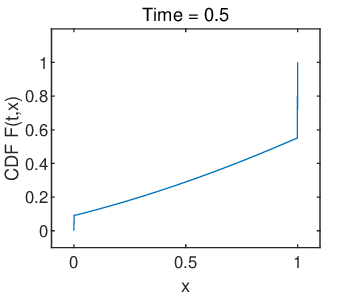}}
{\includegraphics[scale=0.8]{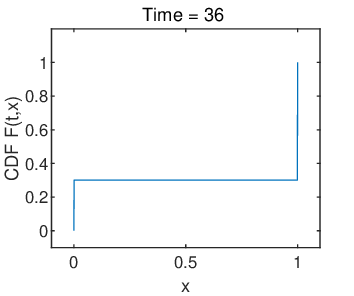}}\\
\caption{Evolution of $F(t,x)$ in {\bf Example 2} with $h=1/1000$, $\tau=1/1000$.}
\label{fig:Ex2Density}
\end{figure}

\begin{figure}[ht]
\centering
\includegraphics[scale=0.8]{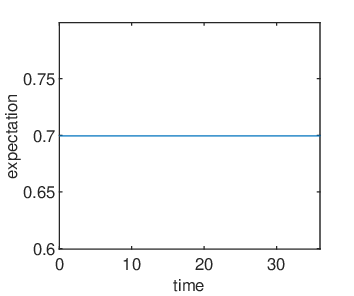}
\caption{Evolution of expectation  in {\bf Example 2} with $h=1/1000$, $\tau=1/1000$.}
\label{fig:Ex2Exp}
\end{figure}

In Fig. \ref{fig:Ex2Density}, the evolution of CDF $F(t,x)$ are shown with partition $h=1/1000$, $\tau=1/1000$. The  discontinuity seems to develop at the boundary  as time evolves. That means Dirac $\delta$ singularities  develop at the boundary points for the original PDF $f$. To verify that the scheme can catch the height of the discontinuous jump, i.e., the fixation  probability, we present the results in Table  \ref{table:EX2}  on different spacial  grid sizes ($h=1/100,
 1/200,1/400, 1/800)$ with a very fine  time step $\tau=1/10000$ at sufficiently large time $T=36$. It can be found that the discontinuity occurs at boundary points and   the  height of the jump on  the two ends   approach  to  the fixation  probability given in Theorem \eqref{thm} with a rate of  $1$st order. In Fig. \ref{fig:Ex2Exp},  the discrete expectation in \eqref{exp-dis} keeps conservation as time evolves  with $h=1/1000$, $\tau=1/1000$.
%
%
\begin{table}[ht]
\centering
\footnotesize{\caption{Discontinuity at the boundary points at  $T=36$  with $\tau=1/10000$ in {\bf Example 2}.}\label{table:EX2}
\begin{threeparttable}
\begin{tabular}{ccccccc}
\hline

\hline
 $h$     &$F_1^N-F_0^N$ &$e_{left}$&Order   & $F_K^N-F_{K-1}^N$&$e_{right}$& Rate \\\hline
 1/100 &0.303030&3.03030e-03&        &0.696970&3.03030e-03&\\\hline
 1/200  &0.301508&1.50754e-03&1.00727&0.698492&1.50753e-03&1.00727\\\hline
 1/400 &0.300752&7.51880e-04&1.00362&0.699248&7.51880e-04&1.00362 \\\hline
 1/800&0.300375&3.75469e-04&1.00181&0.699624&3.75469e-04&1.00180\\
 \hline
\end{tabular}
\begin{tablenotes}
\footnotesize
\item[1]$a^{\infty}=1-x_0=0.3$, $b^{\infty}=x_0=0.7$ in \eqref{errorl}-\eqref{errorr}.
\end{tablenotes}
\end{threeparttable}}
\end{table}

\noindent{\bf Example 3. Natural selection}

In this example, we consider the natural selection case: $$M(x)=x(1-x)(\eta x+\beta).$$ Let $$\mathcal{E}^{\theta}(t):=\int_0^1f(t,x)\theta(x)dx.$$
The  conserved quantity \eqref{equ:Conservation2} is $\mathcal{E}^{\theta}(t) =\mathcal{E}^{\theta}(0)$ with
 {\bf $$\theta(x)= \frac{\int_0^x e^{-\frac{\eta}{2}\bar{x}^2-\beta \bar{x}}d\bar{x}}{\int_0^1e^{-\frac{\eta}{2}x^2-\beta x}dx},$$ } being the solution of \eqref{equ:theta}.

By the integration by parts,  we have
 \begin{equation}
 \begin{aligned}
\mathcal{E}^{\theta}(t)
&=\int_0^1\partial_x F(t,x)\theta(x) dx=F(t,1)\theta(1)-F(t,0)\theta(0)-\int_0^1  F(t,x)\theta'dx\\
&=1-\int_0^1  F(t,x)\frac{e^{-\frac{\eta}{2}x^2-\beta x}}{\int_0^1e^{-\frac{\eta}{2}x^2-\beta x}dx} dx.
\end{aligned}
\end{equation}
Then a discrete conserved quantity can be defined as
$$\mathcal{E}^{\theta,n}_h=1-\frac{1}{A}\left(h\sum\limits_{i=1}^{N-1}F_i^n e^{-\frac{\eta}{2}x_i^2-\beta x_i}+\frac{h}{2}F_0^n e^{-\frac{\eta}{2}x_0^2-\beta x_0}+\frac{h}{2}F_{N}^n e^{-\frac{\eta}{2}x_N^2-\beta x_N}\right),$$
 where $A$ is an approximation of $\int_0^1e^{-\frac{\eta}{2}x^2-\beta x}dx$,
$$A:=h\sum\limits_{i=1}^{N-1} e^{-\frac{\eta}{2}x_i^2-\beta x_i}+\frac{h}{2} e^{-\frac{\eta}{2}x_0^2-\beta x_0}+\frac{h}{2} e^{-\frac{\eta}{2}x_N^2-\beta x_N}.$$

\begin{figure}[ht]
\centering
{\includegraphics[scale=0.8]{PureDensity_t=0.eps}}
{\includegraphics[scale=0.8]{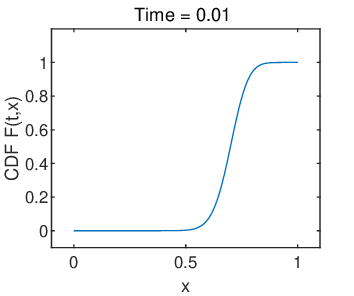}}\\
{\includegraphics[scale=0.8]{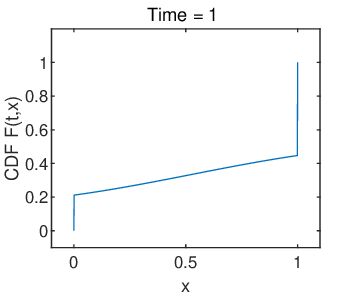}}
{\includegraphics[scale=0.8]{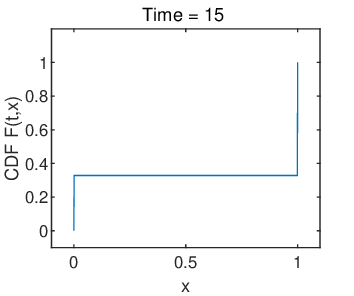}}\\
\caption{Evolution of the numerical solution of CDF in {\bf Example 3} with $h=1/1000$, $\tau=1/1000$.}
\label{fig:Ex3Density}
\end{figure}

\begin{figure}[htbp]
\begin{minipage}[t]{0.5\linewidth}
\centering
\includegraphics[scale=0.8]{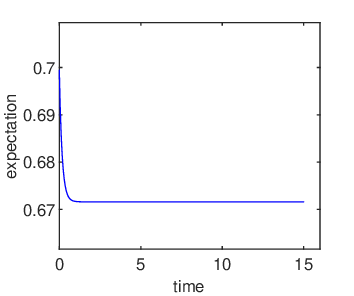}
\caption{Evolution of expectation }\label{fig:Exa3Exp}
\end{minipage}
\begin{minipage}[t]{0.5\linewidth}
\centering
\includegraphics[scale=0.8]{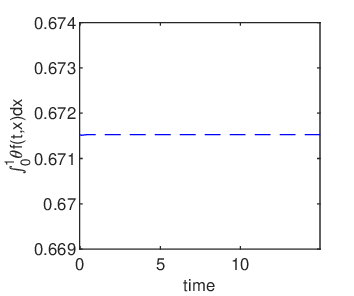}
\caption{Evolution of $\int_0^1 \theta(x)f(t,x)dx$ }\label{fig:Exa3Moment}
\end{minipage}
\end{figure}

In this example,  $\eta=-4$ and $\beta=2$
and  the initial state is given in  \eqref{equ:Nini} with $x_0=0.7$ and $\sigma=0.01$.
 Fig. \ref{fig:Ex3Density} shows that the discontinuous points  emerge at the boundary points.  As time evolves,  the discrete expectation $\mathcal{E}_h^n$ does not keep the conservation and tends to a certain value ($\approx 0.671595$), but  $\mathcal{E}^{\theta,n}_h$ keeps the conservation  and its value is about $0.671529$ shown in Fig. \ref{fig:Exa3Exp}.
Table \ref{Ex3table:boun} shows the ability to catch the jump of the discontinuity  at boundary points under different spacial  grid sizes ($h=1/100,
 1/200,1/400,1/800$), and a  very fine time step $\tau=1/10000$  at $T=15$  when the steady state is approaching.
 The fixation  probability is predicted in a $1$-order accuracy  at left and right boundary point.
%
%


\begin{table}[ht]
\centering
\footnotesize{\caption{Discontinuity at the boundary points at  $T=15$  with $\tau=1/10000$ in {\bf Example 3}.}\label{Ex3table:boun}
\begin{threeparttable}
\begin{tabular}{cccccccc}
\hline

\hline
$h$     &$F_1^N-F_0^N$  &$e_{left}$&Rate & $F_K^N-F_{K-1}^N$ &$e_{right}$&Rate \\\hline
 1/100 &0.330230&2.16298e-03 &  &0.669770&2.16298e-03 & \\\hline
 1/200 &0.329103&1.03683e-03&1.06085 &0.670896&1.03683e-03& 1.06085\\\hline
 1/400 &0.328556 &4.89322e-04&1.08332  &0.671444 &4.89322e-04&1.08332\\\hline
 1/800 &0.328287&2.20144e-04& 1.15234 & 0.671713&  2.20144e-04&1.15234\\
 \hline
\end{tabular}
\begin{tablenotes}
\footnotesize
\item[1]
$e_{left}$ and $e_{right}$ are defined in \eqref{errorl} and \eqref{errorr},  with $b^{\infty}=\mathcal{E}^{\theta,0}_{h^*}$ and $a^{\infty}=1-b^{\infty}$, where the conserved quantity $\mathcal{E}^{\theta,0}_{h^*}=0.671933$  under a very fine space step $h=1/10000$.
\end{tablenotes}
\end{threeparttable}}
\end{table}

\subsection{Mutation case}\label{sec:numericalResults3}
In this section, we discuss the mutation case $M(X)=\gamma(1-x)-\mu x$, $\mu,\gamma\ge 0$, including one-way mutation and two-way mutation.

\noindent {\bf Example 4.  One-way mutation: Muller's ratchet}

\begin{figure}[ht]
\centering
{\includegraphics[scale=0.8]{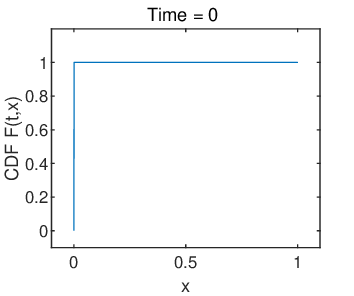}}
{\includegraphics[scale=0.8]{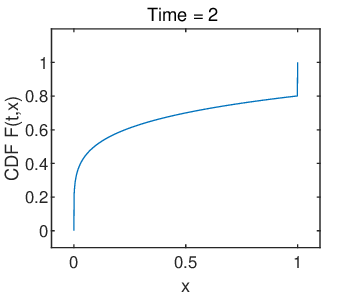}}\\
{\includegraphics[scale=0.8]{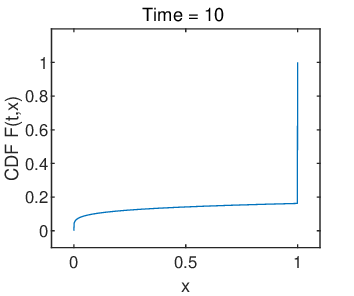}}
{\includegraphics[scale=0.8]{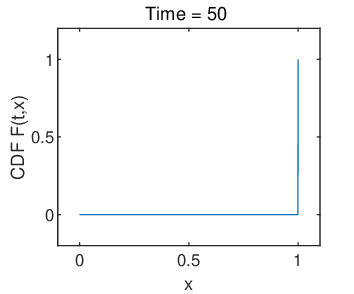}}\\
\caption{Evolution of $F(t,x)$ in {\bf Example 4}  for $\gamma=0.2,\  \mu=0$ with $h=1/1000$, $\tau=1/1000$.}
\label{fig:Ex4DensityMu}
\end{figure}

 One-way mutation  $M(x)=\gamma(1-x)$ with  $\gamma=0.2$ is considered in this example.
The  initial state is taken as
\begin{equation}\label{mutationIni1}
f(0,x)=
\left\{
\begin{aligned}
&\delta_0,\ x=0,\\
&0,\ x\in(0,1],
\end{aligned}
\right.
\end{equation}
i.e., there is a point measure with the whole probability $1$  at $x=0$  at the initial time. That means only the fittest gene $B$ exists in the system.
 Fig. \ref{fig:Ex4DensityMu} shows the evolution of CDF  $F(t,x)$. The   discontinuity  develops at   $x=1$ and the height of the discontinuous jump  rises up to  $1$ eventually.
Table \ref{table:Ex4} shows  the numerical results  with different spatial  step sizes ($h=1/100,\ 1/200,\ 1/400, \ 1/800$) and  a very fine time step $\tau=1/10000$ at  $T=50$, when the steady state is approaching.  It can be found that  the  discontinuity  only emerges at $x=1$ with height of $1$ and no discontinuity happens at $x=0$.
This fact accords with the  Muller's ratchet theory: all fittest gene $B$ will mutate irreversibly to the deleterious gene $A$.

\begin{table}[ht]
\centering
\caption{Discontinuity  at the right boundary point   at $T=50$  in {\bf Example 4}.}\label{table:Ex4}
\footnotesize{\begin{threeparttable}\begin{tabular}{ccc}
\hline

\hline
$h$    &$F_1^N-F_0^N$    & $F_K^N-F_{K-1}^N$ \\\hline
 1/100 &2.22716e-05&0.999946 \\\hline
 1/200 &1.93788e-05&0.999946\\\hline
 1/400  &1.68642e-05&0.999946\\\hline
  1/800& 1.46777e-05&0.999946\\
 \hline
\end{tabular}
\end{threeparttable}}
\end{table}
%

%
\noindent{\bf Example 5. Two-way mutation}

\begin{figure}[ht]
 \centering
{\includegraphics[scale=0.8]{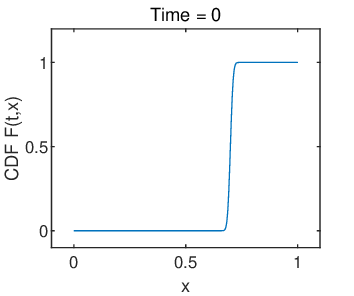}}
{\includegraphics[scale=0.8]{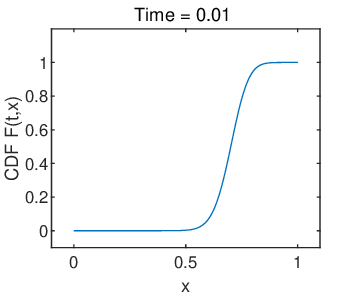}}\\
{\includegraphics[scale=0.8]{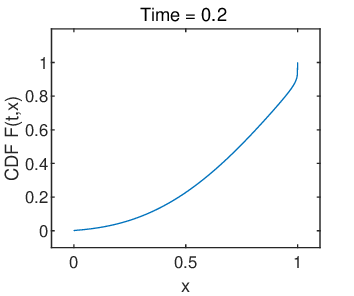}}
{\includegraphics[scale=0.8]{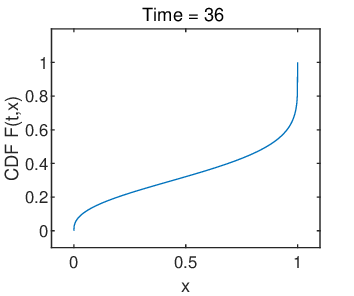}}
\caption{Evolution of $F(t,x)$ in {\bf Example 5} for $\mu=0.2,\ \gamma=0.4$ with $h=1/1000$, $\tau=1/1000$.}
\label{fig:Ex5DensityMu}
\end{figure}

\begin{figure}[ht]
\centering
\includegraphics[scale=0.8]{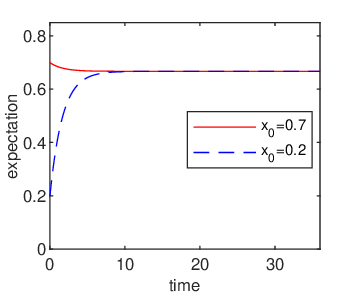}
\caption{Evolution of expectation  in {\bf Example 5} with $\mu=0.2,\ \gamma=0.4$ under $h=1/1000$, $\tau=1/1000$.}
\label{fig:Ex5Exp}
\end{figure}

\begin{figure}[ht]
 \centering
{\includegraphics[scale=0.8]{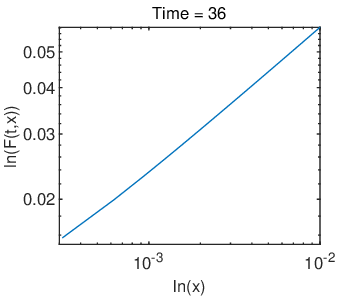}}
{\includegraphics[scale=0.8]{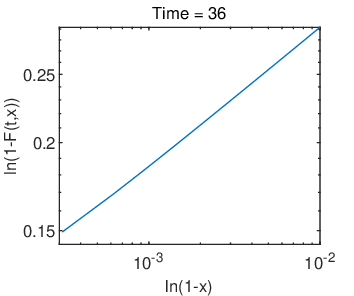}}
\caption{ $\ln(F(t,x))$ near boundary points in {\bf Example 5} at $T=36$ with $h=1/3600$, $\tau=1/10000$.}
\label{fig:Ex5MuPower}
\end{figure}

In this example, we consider a two-way mutation  $M(X)=\gamma(1-x)-\mu x$ with $\mu=0.2,\ \gamma=0.4$. Taking into account that  $F(t,x)$ may be smooth on $x\in [0,1]$, the probability  density $f(t,x)$ can be recovered by \eqref{equf} except $i=1,K-1$.   Central difference is available now,
\begin{equation}\label{equf_TwoWay}
f_i^n=
\left\{
 \begin{aligned}
 &  \frac{F_{i+1}^n-F_{i-1}^n}{2h},\  i=1,\\
 & \frac{F_{i+1}^n-F_{i-1}^n}{2h}, \ i=K-1.\\
 \end{aligned}
 \right.\
\end{equation}
The initial function is chosen as \eqref{equ:Nini} with $\sigma=0.01$.  Fig. \ref{fig:Ex5DensityMu} shows that  the  CDF  may be  continuous as time evolves with $x_0=0.7$  under $h=1/1000$, $\tau=1/1000$.  Fig. \ref{fig:Ex5Exp} shows the expectation is not   conserved as time evolves and towards   to the same value $0.66677$ with different initial $x_0=0.7, 0.2$.   The  left   figure of Fig. \ref{fig:Ex5MuPower} shows  the relationship between $\ln(F(t,x))$ and $\ln(x)$ is approximately linear at $T=36$ with $h=1/3200$ and $\tau=1/10000$.   $\ln(1-F(t,x))$ and $\ln(1-x)$ are also approximately linear in the right figure of Fig. \ref{fig:Ex5MuPower}. The results show that $F(t,x)$ can be  approached by polynomial with $x^{\xi}$ near $x=0$ and $(1-x)^{\eta}$  near $x=1$, where $\xi,\eta$ are positive constant, at $T=36$ near the steady state.

Table \ref{table:Ex5} shows  the behavior of the power law near the boundary points $x=0,1$  with  different initial states ($x_0= 0.7$ and $x_0=0.2$) under different  spatial mesh steps  and the refine time step $\tau=1/10000$ at time $T=36 $. The value of  $F_1^N-F_0^N$ means the point measure does not emerge at the boundary point. The results also show that
the numerical  $F(t,x)$ can be  approached by polynomial $x^{\gamma}$ with  $\gamma\approx 0.4$  near $x=0$ and $(1-x)^{\mu}$  with $\mu\approx 0.2$  near $x=1$, respectively. That means the corresponding probability density $f(t,x)$       is close to $x^{-0.6}$ at $x=0$ and $(1-x)^{-0.8}$  at $x=1$.
It suggests that the numerical results are consistent with theoretical results \eqref{MutationS}.  In addition,  numerical results also show the fact that the steady state has nothing with  the   initial states.
\begin{table}[ht]
\centering
\footnotesize{
\caption{Behavior of power law at boundary points at   $T=36$  with $\tau=1/10000$ in {\bf Example 5}.}\label{table:Ex5}
\begin{threeparttable}
\begin{tabular}{ccccccc}
\hline
 &$x_0=0.7$& &&    \\\hline
$h$    &$F_0^N$ &$F_1^N$  &$\gamma$& $F_{K-1}^N$ &$F_K^N$ &$\mu$ \\\hline
 1/200 &0.00000  &4.78993e-02&                             &7.37404e-01 &1.00000 &\\\hline
 1/400  &0.00000  &3.62271e-02& 0.402935                  &  7.72496e-01&1.00000&   0.206953 \\\hline
 1/800  &0.00000  &2.74198e-02& 0.401851                   &8.02473e-01&1.00000&   0.203842 \\\hline
 1/1600 &0.00000  &2.07643e-02& 0.401110                     &8.28292e-01&1.00000& 0.202093\\\hline
 1/3200 &0.00000  &1.57294e-02& 0.400639                &8.50637e-01&1.00000& 0.201136\\\hline
&$x_0=0.2$& && \\\hline
$h$    &$F_0^N$ &$F_1^N$  &$\gamma$& $F_{K-1}^N$ &$F_K^N$ &$\mu$ \\\hline
 1/200 &0.00000   &4.78993e-02&                             &7.37404e-01 &1.00000 &\\\hline
 1/400  &0.00000  &3.62271e-02& 0.402935                  &  7.72496e-01&1.00000&   0.206953 \\\hline
 1/800  &0.00000  &2.74198e-02& 0.401851                   &8.02473e-01&1.00000&   0.203842 \\\hline
 1/1600 &0.00000  &2.07643e-02& 0.401110                     &8.28292e-01&1.00000& 0.202093\\\hline
 1/3200 &0.00000  &1.57294e-02& 0.400639                &8.50637e-01&1.00000& 0.201136\\\hline
\end{tabular}
\begin{tablenotes}
\footnotesize
\item[1]  $\gamma$ and $\mu$ are obtained by $$\gamma=\ln(F_1^{N}/F_1^{2N})/\ln(2),$$
$$\mu=\ln\left(\frac{F_K^N-F_1^{N}}{F_K^{2N}-F_1^{2N}}\right)/\ln(2).$$
\end{tablenotes}
\end{threeparttable}}
\end{table}
\section{Some discussions about the revised FDM and the standard FDM}\label{sec:FDM2}
 In this section,   we discuss what happens if the revised treatment \eqref{equ:aboundary} is not introduced.

Recalling that $ a(x) =x (1-x)$, the standard finite difference scheme, referred as {\bf sFDM}, is as follows.
Given $F^{n-1}$,   $F^{n}=(F^{n}_0,\cdots,F^{n}_{K})$ such that
\begin{equation}\label{equ:FDM2}
\frac{F^{n}_i-F^{n-1}_i}{\tau}-\frac{ a(x_{i+\frac{1}{2}})D_h F_{i+1}^{n}- a(x_{i-\frac{1}{2}})D_hF_{i}^{n}}{h}+M(x_i) D^{up}_h F^{n}_i=0,\ i=1,\cdots,K-1,
\end{equation}
 subject to $F^{n}_0=0$, \ $F^{n}_{K}=1$, $n=1,\cdots,N$.

\begin{figure}[ht]
\centering
{\includegraphics[scale=0.8]{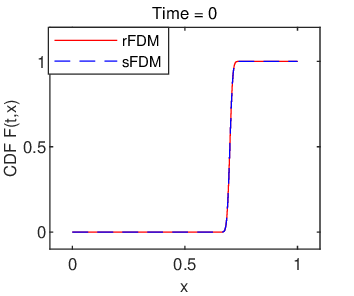}}
{\includegraphics[scale=0.8]{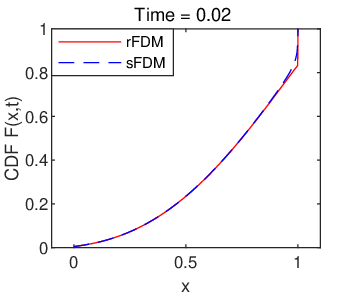}}\\
{\includegraphics[scale=0.8]{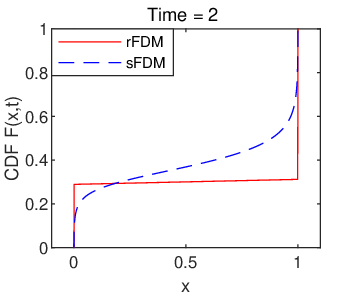}}
{\includegraphics[scale=0.8]{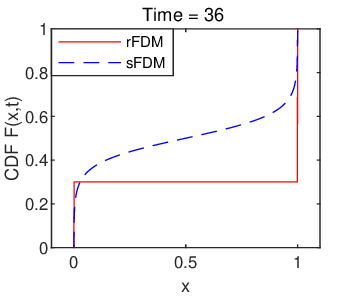}}
\caption{Evolution of $F(t,x) $ by rFDM and sFDM    with $h=1/1000$, $\tau=1/1000$.   }
\label{fig:Ex6DensityFDM1}
\end{figure}

\begin{figure}[ht]
\centering
{\includegraphics[scale=0.8]{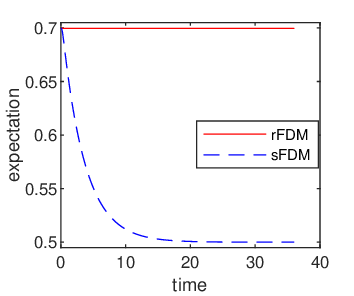}}
\caption{Evolution of expectation in for rFDM and sFDM  with $h=1/1000$, $\tau=1/1000$.}
\label{fig:Ex6Expectation}
\end{figure}

Firstly, we compare  the numerical behavior of the two schemes. Without loss of generality, we consider the pure drift case $M(x)=0$
and take the initial state as \eqref{equ:Nini} with $x_0=0.7$ and $\sigma=0.01$. Numerical results are presented in Figs \ref{fig:Ex6DensityFDM1} and \ref{fig:Ex6Expectation}.

 The evolution of CDF $F(t,x)$ by rFDM \eqref{equ:FDM}-\eqref{equ:aboundary} and sFDM \eqref{equ:FDM2} are shown in Fig. \ref{fig:Ex6DensityFDM1}. As time evolves, the   discontinuity  emerges at the ends $x=0,1$ by rFDM   and the fixation phenomenon is correctly predicted. For sFDM, no evidence implies the  development of the discontinuity, i.e., sFDM fails to predict the fixation phenomenon. To make it more clear, more results by sFDM are presented in Table \ref{table:Ex6} with  different spatial grid size ($h=1/100,\ 1/200,\ 1/400, \ 1/800$) and the fixed time step size $\tau=1/10000$.  It is obvious that sFDM fails to catch the discontinuity that should develop at the ends.

The evolution of expectation in Fig. \ref{fig:Ex6Expectation} shows that   rFDM keeps the conservation of expectation, while sFDM fails.

\begin{table}
\centering
\footnotesize{
\caption{Behavior at boundary points  for   sFDM \eqref{equ:FDM2} at  $T=36$.\label{table:Ex6}}
\begin{threeparttable}
\begin{tabular}{ccccc}
 \hline
 $h$   &$F_1^N-F_0^N$ &$e_{left}$   & $F_K^N-F_{K-1}^N$ &$e_{right}$\\\hline
 1/100 &0.153003&0.146998&0.153003&0.546997 \\\hline
 1/200  &0.138051&0.161949&0.138052&0.561948\\\hline          
 1/400 &0.125864&0.174135&0.125865&0.574135\\\hline             
 1/800&0.115703&0.184297&0.115707&0.584294\\           
 \hline
\end{tabular}
\begin{tablenotes}
\footnotesize
\item[1] $a^{\infty}=1-x_0=0.3$, $b^{\infty}=x_0=0.7$ in \eqref{errorl}-\eqref{errorr}.
\end{tablenotes}
\end{threeparttable}}
\end{table}
The reason why sFDM does not work is that sFDM destroys the rule that the information at boundary points should not be transferred into the domain by diffusion due to the degeneration of the diffusion coefficient $a(0)=a(1)=0$ at the boundary points.

Next, the result of sFDM in Fig. \ref{fig:Ex6DensityFDM1} looks like the one in Fig. \ref{fig:Ex5DensityMu} for two-way mutation case. But what we treat now is the pure drift case $M(x)=0$.  To make this clear, we  take the difference   between   rFDM \eqref{equ:FDM}-\eqref{equ:aboundary}  and sFDM \eqref{equ:FDM2}.
The only difference takes place  at two points $x_1$ and $x_{K-1}$. At $x_1$, sFDM is rFDM plus a term  at the left hand side as
$$ \Lambda_1= a(x_{\frac{1}{2}})D_h F^n_1/h = \frac{1}{2} (1-\frac{h}{2}) D_h F^n_1.$$
This means that at $x_1$, a mutation from gene $B$ to $A$ is numerically introduced to a pure drift case, here $M(x)=\gamma (1-x)$ with a mutation ratio $\gamma=\frac{1}{2}$.
Similarly, at $x_{K-1}$, sFDM is rFDM plus a term at the right hand side as
$$ \Lambda_{K-1}= -a(x_{K-\frac{1}{2}})D_h F^n_{K}/h = -\frac{1}{2} x_{K-\frac{1}{2}} D_h F^n_{K}.$$
That implies that at $x_{K-1}$, a mutation from gene $A$ to $B$ is numerically introduced, here $M(x)=- \mu x$ with a mutation rate $\mu= \frac{1}{2}$.
With these artificial  mutations, fixation can never happen. That's the reason why  sFDM fails to predict the fixation that should happen for pure drift and its numerical results behavior like a problem with two-way mutation.

\section{Conlusions}

We re-model the random genetic drift problem on PDF to a new one on CDF. The possible Dirac $\delta$ singularity on PDF then changes to a discontinuous jump on CDF and the height of the jump is just the fixation probability. The possible singularity of negative power law changes to a bounded positive power law.  A revised finite difference method (rFDM) is proposed to uniformly and effectively handle the pure drift with or without natural selection, one-way mutation and two-way mutation.

The idea working on CDF is a potential way to treat multi-alleles genetic drift problems, where multi-dimensional partial differential equations is involved.
It is quite direct to change the equation on PDF to one on CDF. But it is a challenge now to settle down the boundary condition, which is corresponding to the margin distribution.


\section*{Acknowledgments}
The authors would like to thank  Prof X.F. Chen for very helpful discussions on this topic. C. Duan was supported in part by NSFC 11901109. C. Liu is partially supported by NSF grants DMS-1950868 and DMS2118181. X. Yue was  supported by NSFC 11971342, 12071190 and 12371401.

\bibliographystyle{unsrt}

\end{document}